\documentclass{birkjour}

 \theoremstyle{definition}
 
 \theoremstyle{remark}

 \numberwithin{equation}{section}

\begin{document}

\title[]
{Models from the 19th Century used for\\ Visualizing Optical Phenomena \\and Line Geometry}

\author[]{David E. Rowe}




\maketitle

In 1891, one year after the founding of the Deutsche Mathematiker-Vereinigung (DMV), the Munich mathematician Walther Dyck faced a daunting challenge. Dyck was eager to play an active role in DMV affairs,
just as he had since 1884 as professor at Munich's Technische Hochschule (on his career, see \cite{Hashagen2003}). Thus he agreed, at first happily,  
 to organize an exhibition of mathematical models and instruments for the forthcoming annual meeting of the DMV to be held in Nuremberg. Dyck's plan was exceptionally ambitious, though he encountered many problems along the way. He hoped to obtain support for a truly international exhibition, drawing on work by model-makers in France, Russia, Great Britain, Switzerland, and of course throughout Germany. Most of these countries, however, were only modestly represented in the end, with the notable exception of the British. Dyck had to arrange for finding adequate space in Nuremberg for the exhibition, paying for the transportation costs, and no doubt most time consuming of all, he had to edit single handedly the extensive exhibition catalog in time for the opening. At the last moment, however, a cholera epidemic broke out in Germany that threatened to nullify all of Dyck's plans. When he learned that the DMV meeting had been canceled, he immediately proposed that the next annual gathering take place in 1893 in Munich. There he had both helpful hands as well as ready access to the necessary infrastructure. Moreover, he had already completed the most onerous part of his task: the exhibition catalog
\cite{Dyck1892}.

This paper is concerned with events that long preceded Dyck's exhibition. It does, however, 
deal with some of the models found in one fairly small section of his catalog, entitled line geometry. Under this heading, Dyck subsumed five topics: caustic curves, ray systems, caustic surfaces, evolute surfaces, and line complexes. All except the fourth topic will be discussed below, although with attention to only a few of the models in Dyck's exhibition. My main focus in the final sections of this paper will be on models of co-called complex surfaces. These are special fourth-degree surfaces that Julius Pl\"ucker  introduced in the 1860s for visualizing the local structure of a quadratic line complex. Pl\"ucker's complex surfaces turned out to be closely related to Kummer surfaces \cite{Klein1874} and \cite{Hudson1905}, and both of these types of quartics are examples of caustic surfaces, which arise in  
 geometrical optics. Indeed, 
Kummer surfaces represent a natural generalization of the wave surface, first introduced by Augustin Fresnel  (1788--1827) to explain double refraction in biaxial crystals.

\begin{figure}[ht]
        \centering 
        \includegraphics[width=6.5cm,height=5.0cm]{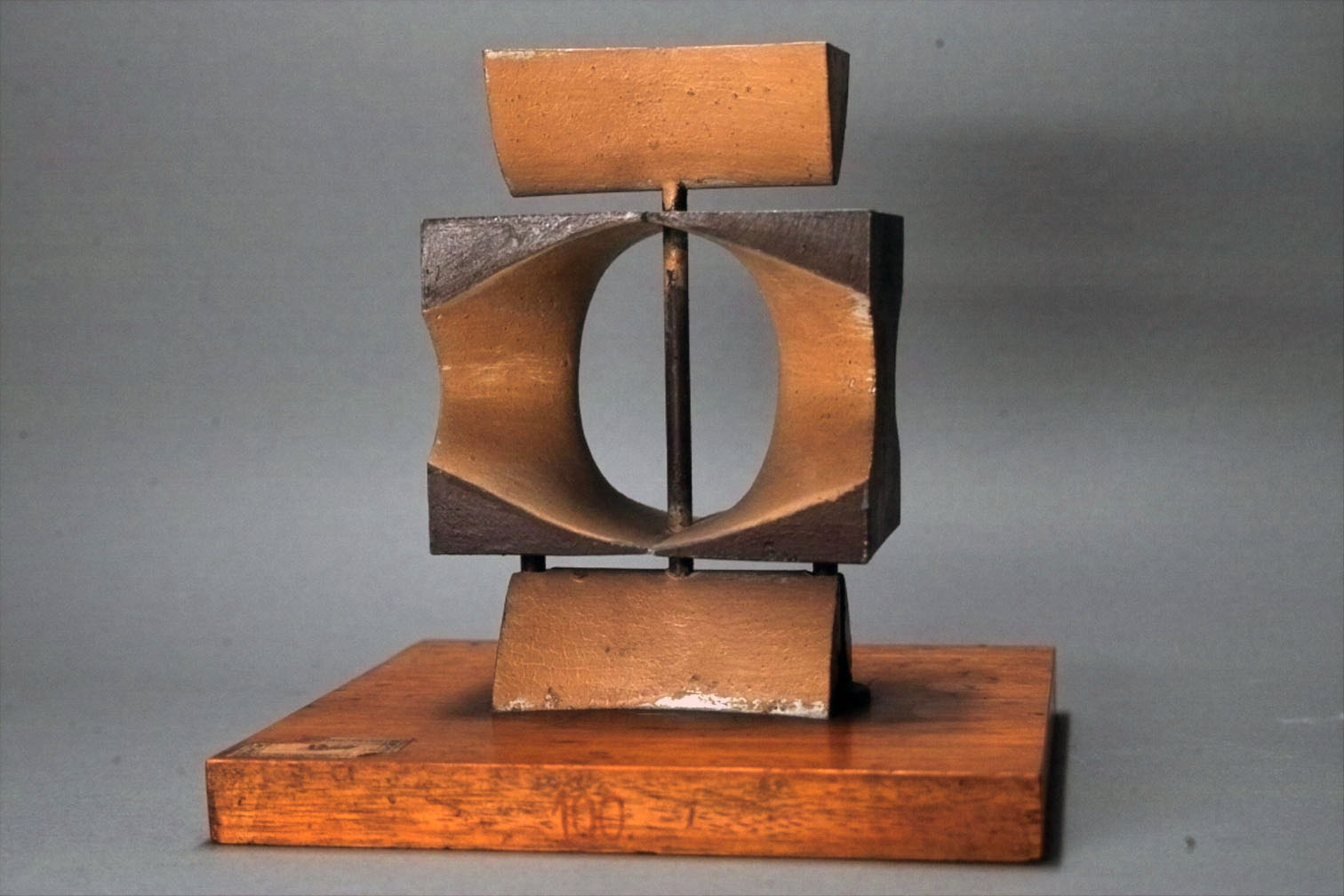}
        \caption{A Pl\"ucker Model from the G\"ottingen Models Collection }
\label{fig:nr23}
\end{figure}

\begin{figure}[ht]
        \centering 
        \includegraphics[width=5.5cm,height=6.5cm]{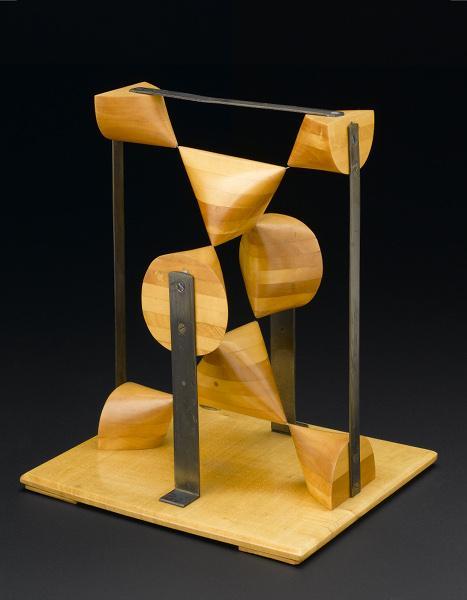}
        \caption{Pl\"ucker's Model of a Complex Surface with 8 Real Nodes, London Mathematical Society}
\label{fig:LMS-13}
\end{figure}

Pl\"ucker  designed  his models to illustrate the various shapes of very special quartic surfaces with a double line (Figure 1 as one example).  These models, made using a heavy metal such as zinc or lead, were displayed for a brief time both in Germany and abroad. In 1866 Pl\"ucker took some with him to Nottingham, where he spoke at a conference about these complex surfaces. Arthur Cayley (1821--1895), Thomas Archer Hirst (1830--1892), Olaus Henrici (1840--1918), and other leading British mathematicians quickly took an interest in them, as they were certainly among the more  
 exotic geometrical objects of their day \cite{Cayley1871}. At Hirst's request, Pl\"ucker later sent a set of boxwood models to the London Mathematical Society (Figure \ref{fig:LMS-13}). These have recently been put  on display at De Morgan House, the headquarters of the London Mathematical Society (see 
https://www.lms.ac.uk/archive/plucker-collection),  though without any information about the mathematics behind them. More extensive collections of Pl\"ucker models can be found in G\"ottingen, 
Munich, and Karlsruhe. Those in G\"ottingen are on display in the collection of models at the Mathematics Institute, which can be viewed online, whereas those at technical universities in Munich, and Karlsruhe are held in storage. These Pl\"ucker models are generally recognizable not only due to their exotic shapes but also because they are far heavier than the more familiar geometric models made in plaster.

Until very recently, Pl\"ucker's models had languished in dusty display cases or in cob-webbed cellar rooms. One has to look quite diligently to find any mention of them in historical sources as well. A very brief description, though, can be found in the final entry on line geometry in Dyck's catalog. This was written by Felix Klein, who during the mid 1860s worked closely with Pl\"ucker at the time he first designed 27 models of complex surfaces. Klein's work on line geometry showed how the Pl\"ucker quartics can be treated as degenerate Kummer surfaces, which afterward were regarded as the surfaces of principle interest for the theory of quadratic line complexes.
Since Klein promoted models and visualization of geometric structures nearly all his life, many mathematicians must have known about Pl\"ucker surfaces in the 1890s. Yet, despite the attention his  models attracted during the final decades of the nineteenth century, particularly among British geometers, by the early twentieth they had nearly disappeared without a trace. What follows can be understood at least in part as an attempt to bring these long lost objects back to life. We begin, though, with some backgroundrelating to the longstanding interplay between optical and mathematical knowledge.

\subsection*{Optics stimulating Mathematics simulating Optics}

Geometrical optics, a field of inquiry nearly as old as mathematics itself, was often closely intertwined with purely mathematical investigations. In his {\it Catoptrica}, Heron of Alexandria derived the law of reflection --- that the angle of incidence equals the angle of reflection --- by invoking a minimal principle, namely, that light beams travel along paths that minimize distance. Heron's main interest, though, seems to have been directed toward achieving a variety of surprising optical effects by means of carefully arranged mirrors. He presumably never reflected on the phenomenon of refraction, however, whereas Claudius Ptolemy performed experiments to determine the angles of refraction for light passing from air to other media, such as water and glass. In the tenth century, the Baghdad mathematician Ibn Sahl  
wrote a commentary on Ptolemy's {\it Optics} as well as an independent study in which he obtained a law of refraction equivalent to Snell's law \cite{Rashed1990}. The latter result 
was found experimentally in Leyden at the beginning of the seventeenth century. Fermat then formulated the principle of least time to account for this refraction law, assuming that the speed of light differs in passing from one medium to another. Leibniz then re-demonstrated Fermat's argument in his famous paper of 1684, in which he introduced the rules for his differential calculus. During the decade that followed, Johann Bernoulli made ingenious use of Fermat's principle to solve the famous brachistochrone problem by treating the earth's gravitational field as a medium with a continously varying ``refraction index''. 

Already in antiquity, it was known that the conic sections had numerous optical properties. Thus, rays emanating from one focal point of an ellipse will be reflected so that they gather in the second focal point. Light rays can also be directed so as to gather in a focal point by means of a paraboloid, which reflects  light rays parallel to its axis into the focus of the surface.
The partially extant work of Diocles, 
{\it On burning mirrors}, shows this genre of interest. This work exerted a strong influence on optics in the Arabic world, particularly on the most important authority on geometrical optics during the medieval period, Ibn al-Haytham. Fermat's contemporary and rival, 
 Descartes undertook numerous investigations of optical phenomena. His study of the 
 focal  properties of lenses led to his discovery of the special quartic curves known today as Cartesian ovals. 

In general, however,  reflected and refracted rays will form caustic curves rather than converging to a point. Thus, a continuous family of plane light rays will be reflected by a smooth curve so as to envelope a second curve, a catacaustic; if refracted, the corresponding envelope is called a diacaustic. Similarly, one speaks of caustic surfaces formed by envelopes of rays that are reflected or refracted by a smooth surface. A glance at the table of contents in L'H\^opital's famous \emph{Analyse des infiniment petits pour l'intelligence des lignes courbes}  from 1696, the first textbook for the differential calculus, shows that the determination of caustics was a major topic of interest already in the seventeenth century. By the nineteenth century, these optical phenomena were being studied as systems of rays, or congruences of lines, in 3-dimensional space. Ray systems are 2-parameter families of lines that undergo reflections and refraction. In 1808 Etienne Malus introduced the notion of normal congruences, which are families of lines orthogonal to a surface \cite{Atzema1993}, 19--37. The theorem of Malus-Dupin then states that a normal congruence that undergoes any number of reflections and refractions on given surfaces, will remain a normal congruence afterward.

\subsection*{Constructing Fresnel's Wave Surface}

\begin{figure}[h]
        \centering 
        \includegraphics[width=6.5cm,height=5.0cm]{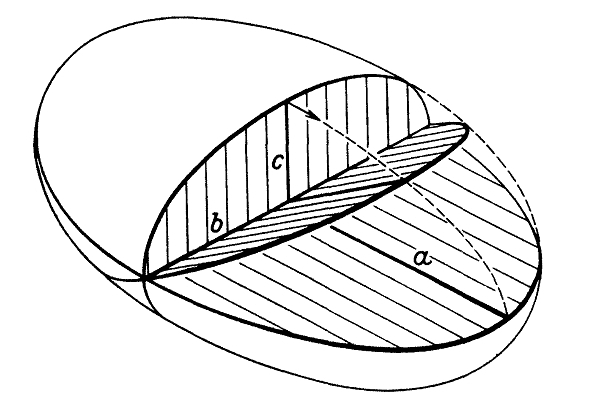}
        \caption{This picture suggests a simple continuity argument that proves why exactly two planes through the center of an ellipsoid will cut the surface in circles. [Hilbert u. Cohn-Vossen 1932], p.  16}
\label{fig:ellipsoid}
\end{figure}

 The Fresnel wave surface 
was an object of considerable interest in the 1820s and long afterward; Pl\"ucker studied it in some detail in \cite{Pluecker1838}.
 Fresnel derived it from a quartic equation for a wave front of light passing through biaxial  crystals, in which the light rays undergo double refraction. This turned out to be a two-leaved surface, though its  deeper geometrical properties eluded him. 
During the 1830s, the Fresnel surface attracted the attention of two natural philosophers at Trinity College Dublin, 
James MacCullagh (1809--1847) and W.R. Hamilton (1805--1865), who revealed several new properties. MacCullagh constructed the wave surface 
directly from Fresnel's index ellipsoid, which represents the strength and orientation of refraction in the crystal \cite{MacCullagh1830}. This 
will be a surface of revolution in the case of  uniaxial crystals, but for biaxial crystals the lengths of the ellipsoid's three axes are unequal. From the classic text \emph{Trait\'e des surfaces du second degr\'e}  by Monge and Hachette (1st ed. 1802),  geometers knew the following two fundamental properties of central quadrics: 1) Each plane section of the surface yields a conic and the conic sections that lie in parallel planes are all similar (i.e., their shape is identical, so they only differ in size). 2) There are two distinguished directions in which the conics cut by parallel planes will be circles. A simple continuity argument shows why property 2) holds (see Figure \ref{fig:ellipsoid}).

\begin{figure}[h]
        \centering 
        \includegraphics[width=6.5cm,height=5.0cm]{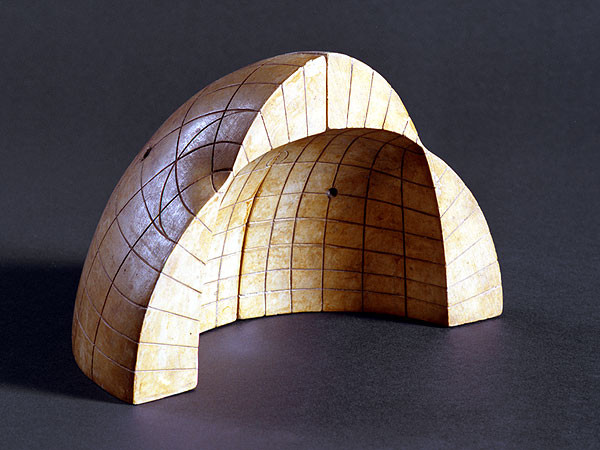}
        \caption{Model of the Interior of a Fresnel wave surface showing two singular points, G\"ottingen Models Collection}
\label{fig:Fresnel}
\end{figure}

\begin{figure}[ht]
\begin{center}
\includegraphics[width=6.5cm,height=5.0cm]{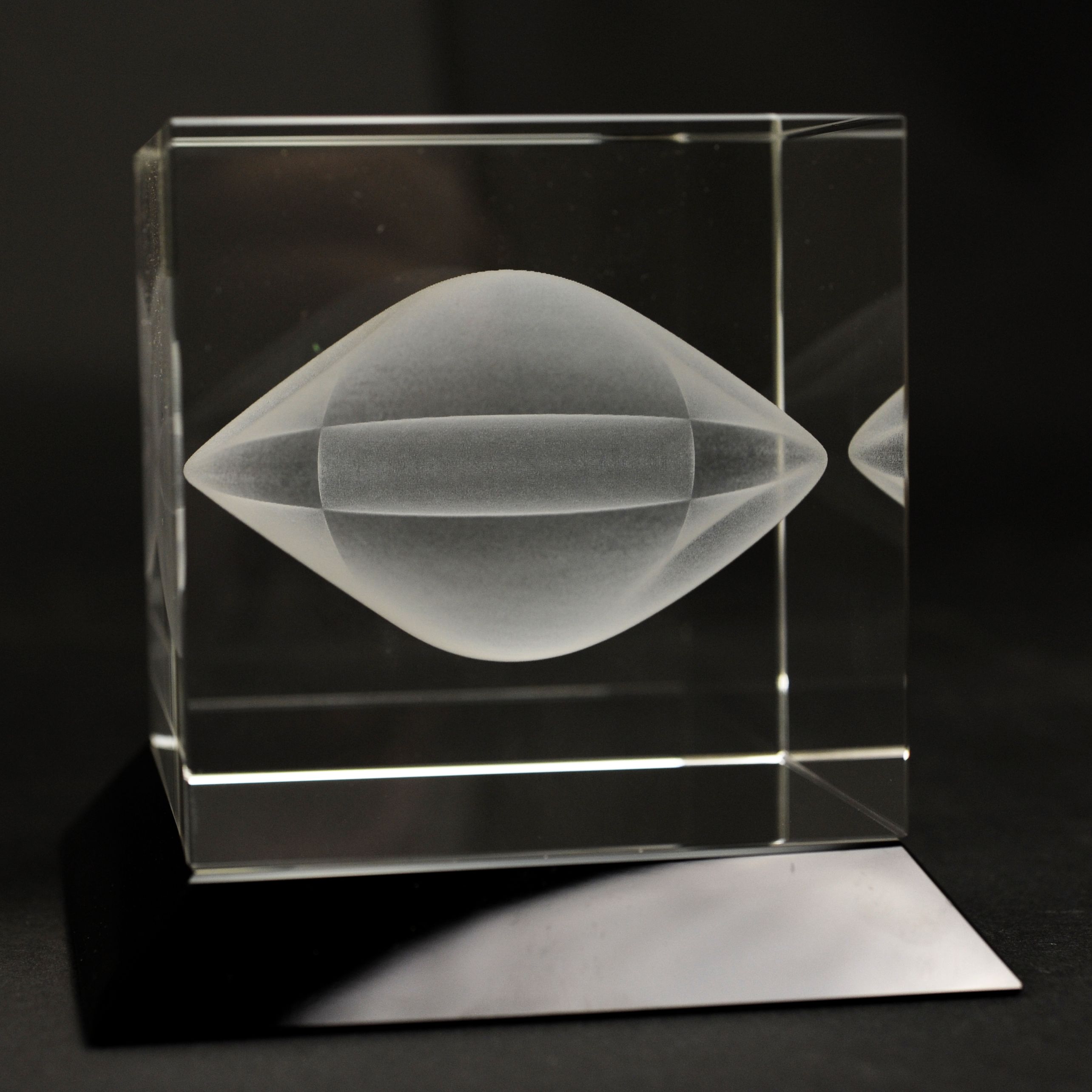}
\caption{Fresnel's Wave Surface as a 3d-laser-in-glass object, produced by Oliver Labs. This recent technique has the advantage over historical models of exposing the surface alone so that its structure, including  
 the four singular points, can be visualized quite adequately. This semi-transparency also makes it possible to see both the outer and inner shells simultaneously; older historical models required two parts in order to be able to look inside.}
\label{fig:Fresnel-Labs}
\end{center}
\end{figure}

These two distinguished directions correspond to the four singular points of the wave surface where its two leaves fall together, as illustrated in Figure \ref{fig:Fresnel}.
 Here the space between the leaves has been filled in, so this model is actually a quarter section of a solid figure whose outer and inner shells belong to the surface. Using laser-in-glass technology, Oliver Labs has produced models of the Fresnel surface that provide a far clearer picture of its geometrical structure. An idea of what such a model looks like can be gathered from Figure \ref{fig:Fresnel-Labs}. This model was designed so as to enhance the  regions surrounding the singular points. The transparent glass makes it far easier to visualize the conical structure formed in passing from a singular point to the outer surface.

MacCullagh's colleague in Dublin, the mathematician-astronomer W.R. Hamilton, was the first to recognize the significance of this geometrical structure for optics. In \cite{Hamilton1837}\footnote{The publication date is misleading, as Hamilton wrote this Third Supplement some years earlier; he discovered conical refraction already in 1832.} he predicted that light rays propagating in two special directions should split into a whole circle of light, a phenomenon that came to be known as conical refraction. Hamilton discussed this theoretical result  in November 1832 with a local expert on experimental optics, Humphrey Lloyd (1800--1881), who confirmed this prediction less than a month later (see \cite{Hankins1980}, 88--95). This shows how a purely mathematical feature of the Fresnel surface led to the discovery of a new physical phenomenon. Hamilton's successful prediction of conical refraction  was celebrated by several contemporary investigators. Pl\"ucker later wrote: ``No experiment of physics ever made such an impression on me \dots it was a thing unheard of and completely without analogy"' (cited from \cite{Hankins1980}, 95). 

MacCullagh, Hamilton, and George Salmon all helped to make the Fresnel wave surface the most prominent example of a quartic surface in the mathematical literature. By introducing complex numbers, Salmon showed that the Fresnel surface actually has 16 singular points, only four of which are real \cite{Salmon1847}. It thus 
belongs to the special class of quartics that would come into prominence with Ernst Eduard Kummer's (1810--1893) work in the 1860s. Salmon  
 showed further that the Fresnel surface was of the fourth order and class with 16 double points and 16 double planes. These turned out to be characteristic properties of the special quartics that Kummer first described in \cite{Kummer1864}.

\subsection*{Constructing Infinitely Thin Pencils of Rays}

Already in \cite{Hamilton1828}, his first study of ray systems, Hamilton dealt with infinitely thin pencils of lines, a topic that would recur throughout much of the nineteenth century \cite{Atzema1993}. He returned to this topic in his Third Supplement \cite{Hamilton1837}, where he described three types of structures that can occur when considering rays of lines infinitely close to a given line. In his earlier work, Hamilton focused his attention entirely on so-called normal systems of rays, those that satisfy the theorem of Malus-Dupin. For such systems, i.e. rays that undergo simple reflection and refraction, only the first type of structure can occur. On a given line line $\ell$, the thin pencil consists of an infinitesimal subset within a given congruence of lines that fill all of space. This pencil determines three special points $M, F_1, F_2$ on $\ell$, the midpoint and two foci, such that $\overline{MF_1}= \overline{MF_2}$. In addition, there are two focal planes, $\Pi_1$ and $\Pi_2$, containing $\ell$, which determine two infinitesimal lines $\epsilon_1, \epsilon_2$ lying, respectively, in $\Pi_1, \Pi_2$, perpendicular to $\ell$, and passing through $F_1$ and  $F_2$. Hamilton then deduced that the infinitely thin pencil centered around $\ell$ can be constructed as the set of lines that intersect both $\epsilon_1$ and $\epsilon_2$. He noted, further,  that in the case of a normal congruence the two focal planes will always be perpendicular to one another. This corresponds to the first type of pencil, whereas for the second type the planes $\Pi_1$ and $\Pi_2$ will no longer be perpendicular. Hamilton's third type of pencil no longer has real foci $F_1, F_2$; instead these are conjugate imaginary points.

\begin{figure}[h]
        \centering 
        \includegraphics[width=9.5cm,height=7.0cm]{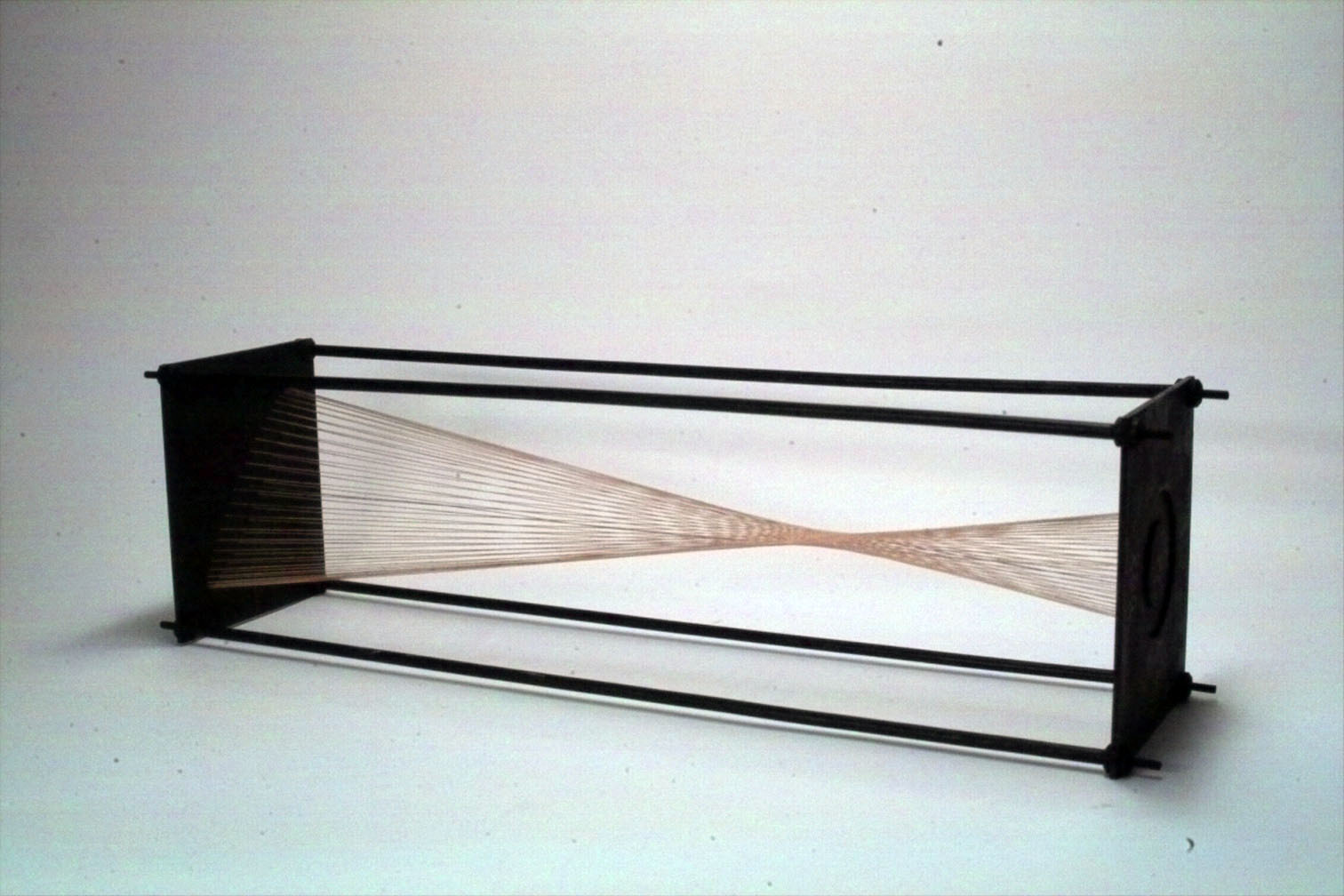}
        \caption{Model by Wilhelm Apel of Hamilton's second type of infinitely thin pencil of lines, G\"ottingen Models Collection}
\label{fig:kp2}
\end{figure}

Few of those who afterward studied infinitely thin pencils of lines took notice of Hamilton's general theory, as pointed out in \cite{Kummer1860a}, which seems to have been the first study to deal with it. One of Kummer's main interests concerned the angle between the two focal 
planes $\Pi_1$ and $\Pi_2$, which would vary within the congruence. He was also interested to consider how the second and third types of pencils arose in uniaxial and biaxial crystals, where the wave fronts form, respectively, ellipsoids of revolution and Fresnel surfaces. Only in the latter case does Hamilton's third type of pencil come into play. In conjunction with this study, Kummer produced three string models for visualizing the structure of each type, describing these in \cite{Kummer1860b}. The models are meant to show the surface of lines that bounds the infinitely thin pencil. Kummer constructed this surface by taking a small circle $\Omega$ in the plane passing through $M$ perpendicular to $\ell$. The lines that meet $\epsilon_1, \epsilon_2$, and $\Omega$ then form a fourth-order surface that bounds the pencil. These models served  as the prototype for others that were produced in G\"ottingen by Wilhelm Apel, whose father had earlier founded a workshop for building scientific instruments in G\"ottingen (see Figure \ref{fig:kp2}).

Apel's three models were among the many on display in Munich in 1893. In the catalog, the Munich geometer Sebastian Finsterwalder briefly explained their construction (\cite{Dyck1892}, 280). Kummer's motivation, as described in \cite{Kummer1860b}, was at least partly physical. He formulated a general theorem, valid for any homogeneous transparent medium, that generalized what he had found for the Fresnel wave surface. In the latter case, the two focal planes are correlated and cut out reciprocal curves which are tiny conics, familiar in differential geometry as the Dupin indicatrix. Kummer asserted that the same type of correlation will occur in general. He also claimed that Hamilton's three types of pencils are the only possible cases, excepting the special types with a conical or cylindrical structure, arising when the lines pass through a point $P\in \ell$, which will be a cylinder if $\overline{PM}= \infty$.

\subsection*{Kummer Surfaces}

Since the Fresnel wave surface arose from a physical problem -- to account for double refraction in biaxial crystals -- its metrical properties guided the initial investigations. Geometers had, however, long recognized that algebraic surfaces can be studied from a purely projective standpoint. This gradually led to the recognition that the Fresnel surface belonged to a large class of quartics, culminating with Kummer's fundamental paper \cite{Kummer1864}. 
In the 1840s 
Arthur Cayley took up the study of special  quartic surfaces he called tetrahedroids \cite{Hudson1905}, 89--92. These have the property that the four planes of a  tetrahedron meet the surface in pairs of conics that pass through four singularities of the surface. The Fresnel surface corresponds to the special case where, in each of the four planes, one of the two conics is a circle. These pairs of curves intersect in 4 points, which are real in only one of the four planes. The other twelve singular points are imaginary, so altogether there are sixteen. In the 1860s, Kummer used elaborate, yet conventional algebraic machinery to explore the properties of 
 quartic surfaces with the maximum possible number of singularities, namely 16 \cite{Kummer1864}. At the same time, he set about producing models of various types of quartics, though, unlike cubics, quartic surfaces proved too plentiful and complicated to succumb to a general classification \cite{Rowe2013}.

\begin{figure}[h]
        \centering 
        \includegraphics[width=8.5cm,height=6.0cm]{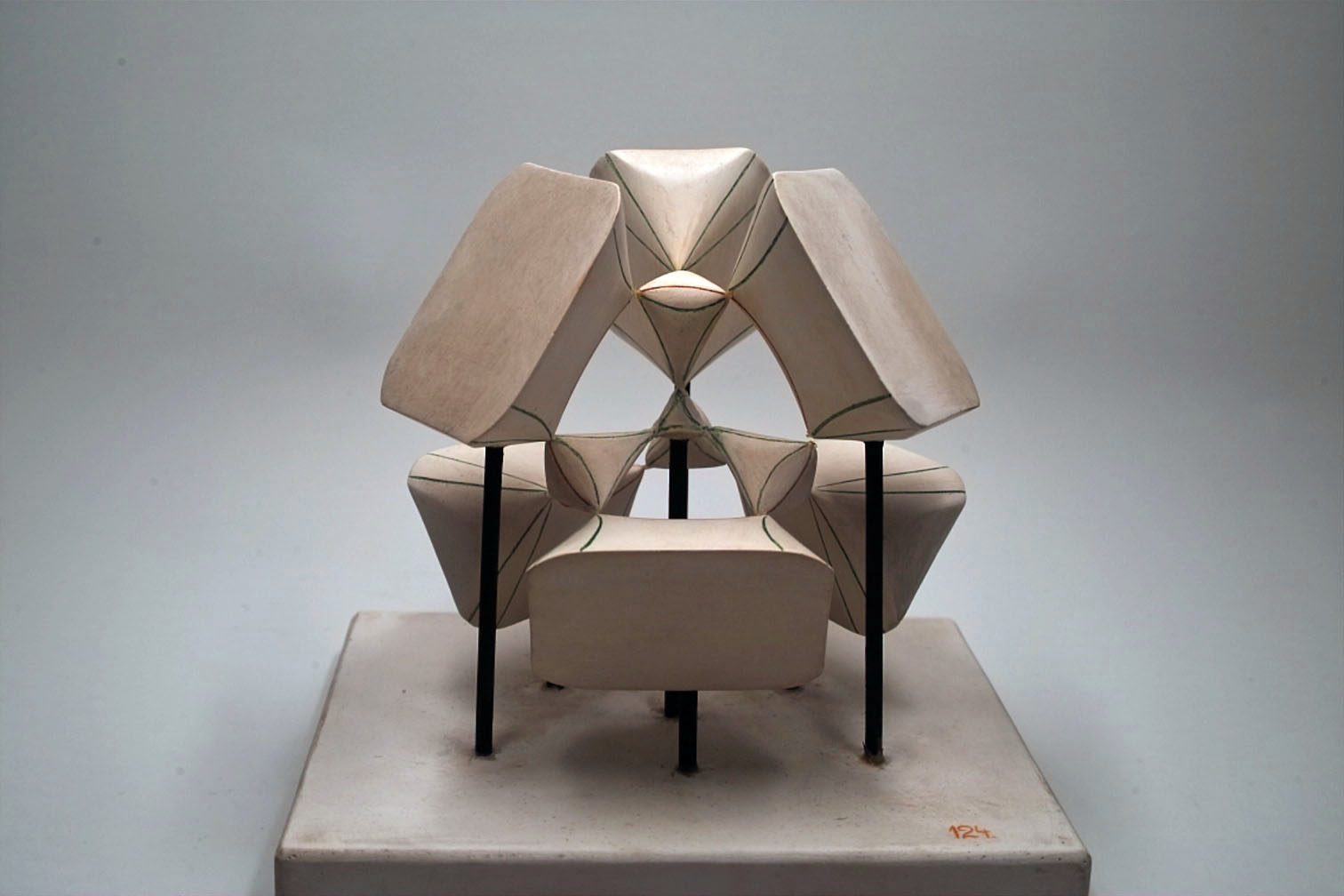}
        \caption{Karl Rohn's model of a Kummer surface with 16 real nodes \cite{Rohn1877}. G\"ottingen Mathematical Institute.}
\label{fig:Rohn_Kummer}
\end{figure}

The 16 singular points of a 
 Kummer quartic form an interesting and much studied spatial configuration. These points  
lie in groups of six in 16 singular planes, or tropes (see Figure 7). Each of these 16 tropes is a tangent plane to the surface that touches it along a conic section which contains six of the 16 singularities. These sets of six points can thus be seen to lie in a special position, since only five coplanar points in general position determine a conic. In fact, the  singular points and planes of a Kummer surface form a symmetric (16,6) configuration, so six singular planes also pass through each of the singular points of the surface. Since all these singularities can be real, this case posed an obvious challenge to model makers, beginning with Kummer himself. The model shown in Figure 7
 was produced in Munich by Klein's student, Karl Rohn; he described it along with two others in \cite{Rohn1877}. Soon afterward, these models proliferated widely, since they were marketed by the Darmstadt firm of Ludwig Brill (Series II, Nr. 1) (for photos, see \cite{Fischer1986}, 34--37). In 1911 these models cost between 21 and 32.50 Marks, which today would be roughly 105 to 160 Euro. 
Brill's brother was Klein's colleague at the TH M\"unchen, Alexander Brill (1842--1935), one of the key figures who promoted model-making during the era. Below, we briefly indicate how the model facilitates the visualization of the (16,6) configuration of singularities on a Kummer surface \cite{Rowe2018b}.

This model represents a finite portion of a highly symmetric Kummer quartic situated in projective 3-space. All of its singularities are visible in the model, which can be seen to consist of eight tetrahedral-like pieces: an inner tetrahedron with four others attached to each vertex, plus three outer tetrahedra that join at infinity. These last three tetrahedra have been split and thus truncated  into two pieces by planar sections, which is why the model has six outer pieces. By imagining the opposite pieces to extend through infinity, they would then join to form the three outer tetrahedra.

 If we now consider the 
 three nearest outer pieces, we see that each contains two vertices that  belong to two of the three tetrahedra associated with the front face of the inner tetrahedron. This gives us six double points, and from the model we can see that these six  points are coplanar. Evidently this is  one of the 16 tropes that touches the surface along a conic, which has been drawn on the model. Similarly, one can  visualize three other tropes by rotating the model to bring into view another set of three half-tetrahedra lying in a plane parallel to one of the three formerly invisible faces of the inner tetrahedron. This procedure produces 4 of the 16 tropes. Notice, furthermore, that the vertices of each of the three outer tetrahedra are shared by one vertex on each of the four middle tetrahedra. This accounts for the full symmetry of the (16,6) incidence configuration of singularities on the  
surface. The remaining 12 tropes can now be brought into view by using projective mappings that interchange the inner tetrahedron with one of the three outer tetrahedra, whereby four new tropes  become visible, thus yielding 16 altogether. This model of a highly symmetric Kummer surface also served as a starting point for Klein's models  of Pl\"ucker's complex surfaces, which he later related to Rohn's important paper \cite{Rohn1879}.

An earlier version of a similar model was built in 1869, when Klein was in Berlin. There he met with his friend Albert Wenker, who designed it for him. Although this Wenker model has disappeared, it should be remembered for the role it played shortly afterward when Klein and Sophus Lie were collaborating for the first time. After meeting during the fall of 1869 in Berlin, where they both attended Kummer's seminar, they continued to work together the following spring in Paris. During this time, Lie found his famous line-to-sphere mapping, a contact transformation with many interesting properties \cite{Lie1872}. Lie used it to map families of lines and the surfaces enveloped by them to families of spheres and their corresponding envelopes. 

\begin{figure}[h]
        \centering 
        \includegraphics[width=6.5cm,height=8.5cm]{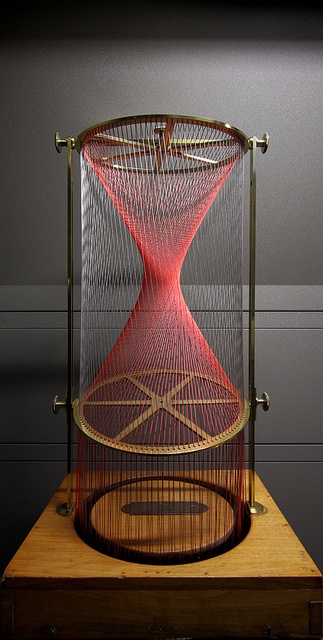}
        \caption{Historical model of a hyperboloid of one sheet by Theodor Olivier.}
\label{fig:Olivier}
\end{figure}

The simplest example arises when the lines are the generators of a hyperboloid of one sheet (Figure \ref{fig:Olivier}), in which case the 
 spheres envelope a Dupin cyclide (Figure 9). One can picture this most easily by taking three skew lines in space, which then map to three spheres. The set of all lines that meet these mutually skew lines form the generators of a quadric surface, and since Lie's mapping is 
a contact transformation these generators go over to a one-parameter family of spheres tangent to three fixed spheres. Since the second system of generators has the same property, the image is a Dupin cyclide, since these arise as surfaces  enveloped by two families of spheres. 
 Moreover, Lie's mapping has a special property of great importance for differential geometry: it sends
the asymptotic curves of one surface to the curvature lines of another. In this simple case, the generators themselves are the asymptotic curves on a hyperboloid, and these are mapped to circles of tangency on the Dupin cyclide (see \cite{Lie/Scheffers1896}, 470--475).

\begin{figure}[h]
        \centering 
        \includegraphics[width=9.5cm,height=6.5cm]{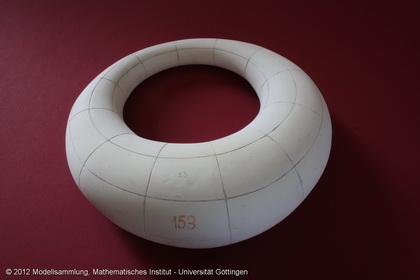}
        \caption{Historical model of a cyclide from the Brill collection.}
\label{fig:cyclide}
\end{figure}

In 1864 Gaston Darboux and Theodore Moutard had begun to work on generalized cyclides, which they studied in the context of inversive geometry (on Darboux's early geometrical work, see \cite{Croizat2016}).
 Klein and Lie learned about this new French theory when they met Darboux just before the outbreak of the Franco-Prussian War. These cyclides are special quartic surfaces with the property that they meet the plane at infinity in a double curve, namely the imaginary circle that lies on all spheres. Darboux also found that their lines of curvature are algebraic curves of degree eight. This finding set up one of Lie's earliest discoveries, communicated to Darboux at that time. This came from Lie's line-to-sphere mapping, when he considered the caustic surface enveloped by lines in a congruence of the second order and class \cite{Rowe1989}.  This was precisely the context that had led 
Kummer to his discovery of Kummer surfaces. Lie found that these Kummer quartics map to generalized cyclides of Darboux, and since the curvature lines of the latter were known, he immediately deduced that the asymptotic curves on a Kummer surface must be algebraic of degree sixteen.

\begin{figure}[h]
        \centering 
        \includegraphics[width=9.5cm,height=6.5cm]{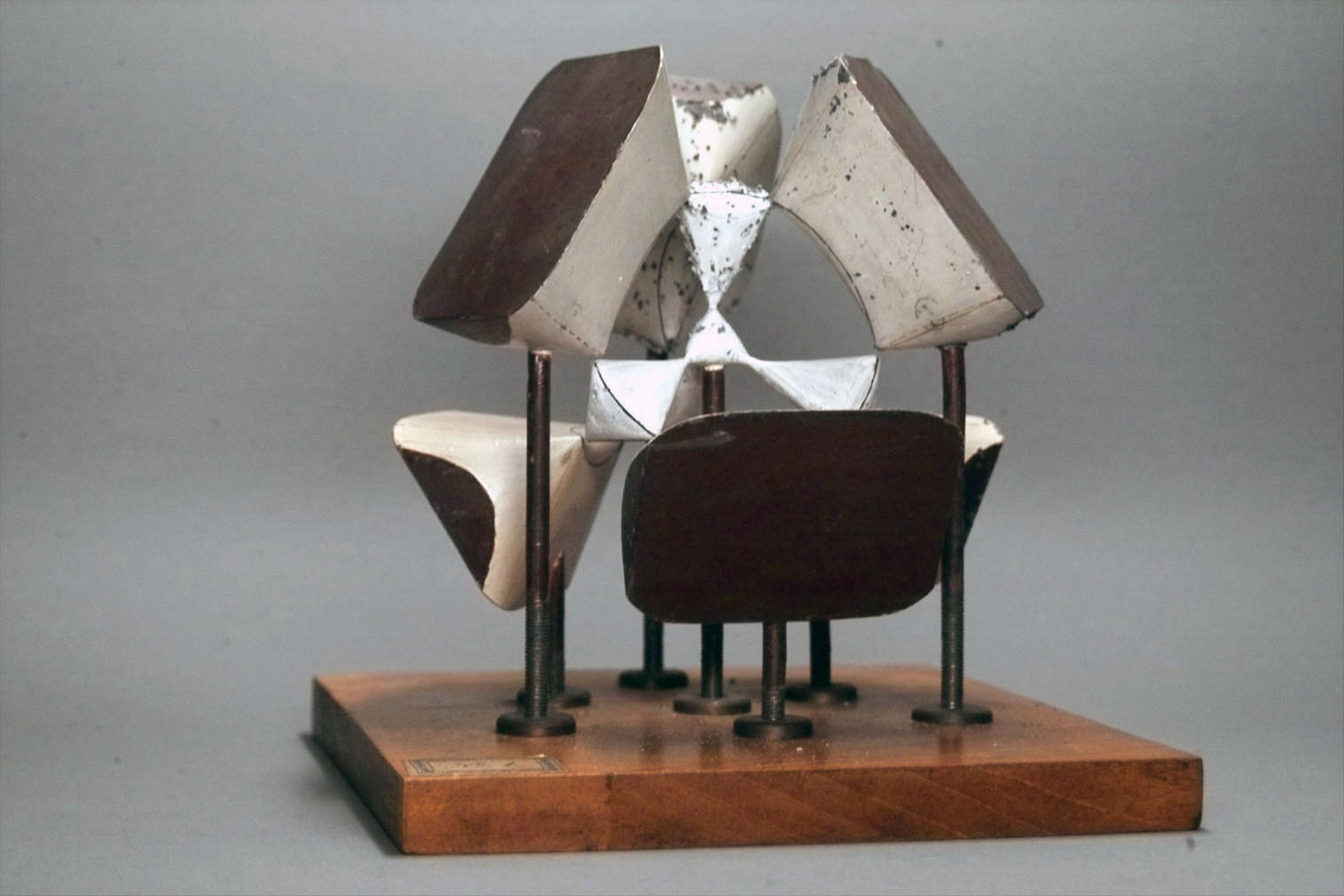}
        \caption{Klein's model of a Kummer surface designed in 1871. Courtesy of the Collection of Mathematical Models,  G\"ottingen University.}
\label{fig:Klein-Kummer}
\end{figure}

\begin{figure}[h]
        \centering 
        \includegraphics[width=6.5cm,height=9.5cm]{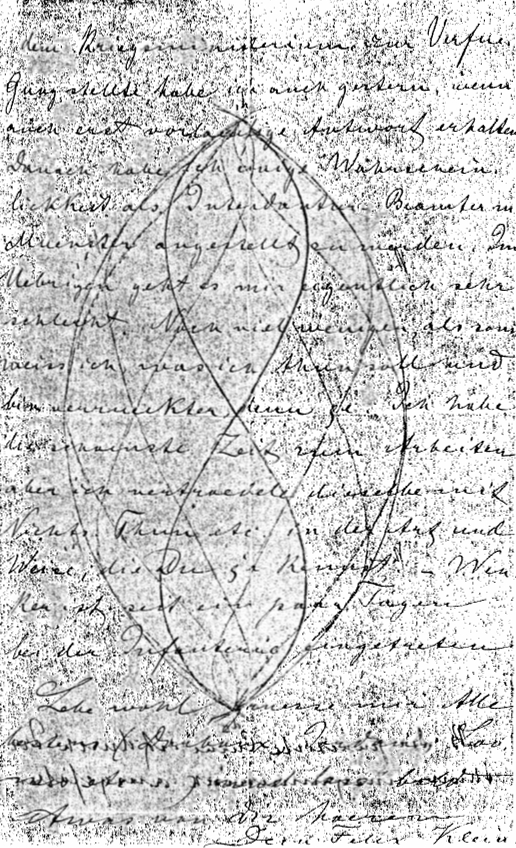}
        \caption{Klein's sketch of the asymptotic curves between two double points on a Kummer surface from his letter to Sophus Lie, 29 July 1870.}
\label{fig:maya}
\end{figure}

Lie communicated these findings to the Norwegian Scientific Society in Christiania in early July 1870, but this note was only published by Ludwig Sylow in 1899, the year of Lie's death (\cite{Lie1934},  86--87). In Paris, Lie and Klein discussed this breakthrough in detail. At first Klein took a skeptical view of Lie's claim,  but then he quickly realized that 
he had already come across these same curves in his own work without recognizing that they were asymptotic curves. 
Little more than a week later, however, Klein had to flee from Paris. Staying at his parents' home in Dusseldorf, he continued to think about the paths of these asymptotic curves as well as their singularities. He did so by tracing them on the physical model  of a Kummer surface made by his friend Albert Wenker. After a short time, 
 Klein realized that what he he had told Lie earlier in Paris about the singularities of these asymptotic curves was, in fact, wrong. After giving the necessary corrections, he wrote:

\begin{quote}

I came across these things by means of Wenker's model, on which I
wanted to sketch asymptotic curves. To give you a sort of intuitive idea
how such curves look, I enclose a sketch. The Kummer surface contains
hyperboloid parts, like those sketched; these are bounded by two of the
six conics ($K_1$ and $K_2$) and extend from one double point ($d_1$) to another
($d_2$). Two of the curves are drawn more boldly; these are the two that not
only belong to linear complexes but also are curves with four-point contact.
They pass through $d_1$  and $d_2$  readily, whereas the remaining curves have
cusps there. This is also evident from the model. At the same time, one sees
how $K_1$  and $K_1$  are true enveloping curves. (Klein to Lie, 29 July 1870)
\end{quote}

By the ``hyperboloid parts'' on a Kummer surface, Klein meant those places where the curvature was negative. Only in these regions are the asymptotic curves real and hence visible. Klein later reproduced the  same figure in the note that he and Lie sent to Kummer for publication in the {\it Monatsberichte} of the Prussian Academy \cite{Klein/Lie1870}. Afterward, this picture became a standard part of the growing literature on Kummer surfaces, as for example in \cite{Rohn1879}, 340 and \cite{Hudson1905}, 119. Klein had a special fondness for it, too. When five years later he married Anna Hegel, granddaughter of the famous philosopher, he ordered a wedding dress decorated with these arabesque curves. No doubt he was happy to explain their significance to any would-be listener.

\subsection*{Pl\"ucker's Complex Surfaces}

During the last years of his life, 
Pl\"ucker  devoted a great deal of attention to studying and building the special quartics he called {\it complex surfaces}. Their name comes from line geometry because these quartics are enveloped by subsets of lines in a quadratic line complex. These are the central objects of 
interest in a field of inquiry he largely invented: line geometry \cite{Pluecker1868}. In classical studies in geometry, the fundamental elements are points, whereas in projective geometry the principle of duality puts points and planes on an equal footing. Pl\"ucker's line geometry, on the other hand, took the manifold of all  
lines in space as fundamental. These, however, require four parameters rather than only three, as with points and planes. One can introduce different systems of line coordinates for this purpose, but in any case the space of all lines will be four-dimensional. In order to 
appreciate Pl\"ucker's motivation for studying complex surfaces, a few basic concepts from 
classical line geometry are needed.

Using line coordinates, an algebraic equation corresponds to a 3-param\-eter family of lines, a {\it line complex}. For degree 2, a simple, yet instructive case comes by taking the lines  
tangent to a nonsingular quadric surface $F_2$, for example an ellipsoid. Here the local structure is immediately obvious: for a given point $P\notin F_2$, the lines through $P$ tangent to $F_2$ form a quadratic cone (real or imaginary). This cone, however, collapses to a tangent plane whenever 
$P\in F_2$, which means that $F_2$ is the associated {\it singularity surface} $S$ of the complex. Dually, a typical plane $\pi$ will cut out a conic, 
$\pi \cap F_2= C_2$, whose tangents are also tangents to $F_2$. The only exceptions are planes $\pi =T_P$ tangent at $P\in F_2$, where the lines enveloping the conic then collapse to a pencil of lines centered at $P$. This shows $S$ is self-dual, since the same surface arises for points as for planes, a property that holds in all quadratic complexes. 

The special character of this elementary quadratic line complex can be seen from its local structure. In the general case, the cone of lines 
 associated with a point $P\in S$ will degenerate into two planar pencils of lines, one centered at $P$ and another at a point $Q$, where the line $PQ$ is a double line belonging to both pencils. Similarly, a singular plane $T_P$ is one in which the conic degenerates into two point pencils, one at $P$ and another at a point $Q$. The tangent cones to a quadric surface $F_2$ thus correspond to the case where the points $P$ and $Q$ fall together forming a single pencil of lines in the plane $T_P$, counted twice. Similarly the 
singularity surface $S$ is the degenerate quartic obtained by doubling $F_2$, so $F_2^2$. In the general case $S$ is a Kummer quartic, which depends on 18 parameters, not just 9, as in the case of an 
an $F_2$. Moreover, the tangents to an $F_2$ completely determine the complex, so it, too, has only 9 parameters. A generic quadratic complex has 19 because, as
 Klein  showed in \cite{Klein1870}, any given Kummer surface is the singularity surface for a one-parameter family of quadratic complexes. 
The same holds for another special case, the 
tetrahedral complexes $T(\lambda)$, which depend on 13-parameters. 
Given any 4 planes in general position, $T(\lambda)$ is the complex determined by those lines whose 4 intersection points with this tetrahedron have a fixed cross ratio $\lambda$. The singularity surface $S$ of $T(\lambda)$  is the tetrahedron itself, and by varying $\lambda$ one gets a one-parameter family of quadratic complexes. 

\begin{figure}[ht]
        \centering 
        \includegraphics[width= 0.33 \textwidth]{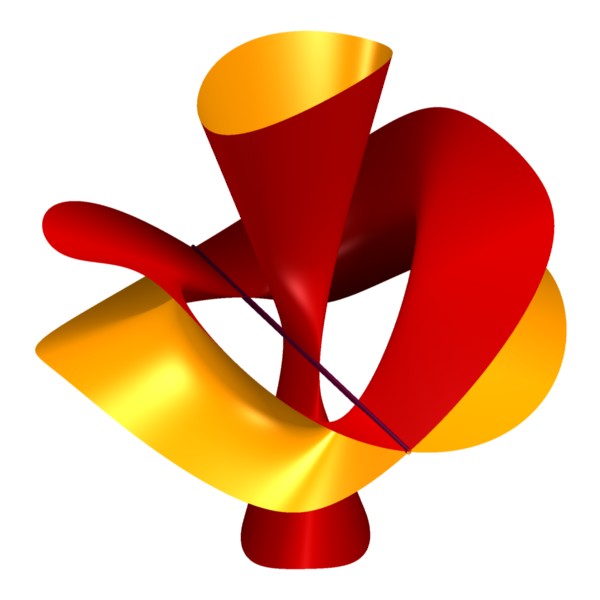}
        \caption{A Pl\"ucker complex surface showing four pinch points on its double line. Graphic by Oliver Labs.}
\label{fig:Pluecker_11112}
\end{figure}

The totality of lines in a quadratic complex is nearly impossible to  visualize, except in special cases like the tangents to a quadric surface or those that belong to a tetrahedral complex. 
Since Pl\"ucker's work was strongly guided by {\it Anschauung} \cite{Clebsch1872}, 
 he studied the local structure of quadratic complexes from the standpoint of lines, the fundamental objects in line geometry \cite{Pluecker1869}.
 His idea was to view the lines in 
a complex $K_2$  with respect to a fixed  line $g$, where in most cases $g\notin K_2$. More precisely, he considered the lines in $K_2  \cap K_1(g)$, where $ K_1(g)$ is the first-degree complex consisting of those lines in space that meet $g$. The intersection of these two complexes produces a 2-parameter family of lines, a line congruence of the second order and class.
This type of ray system was familiar from geometrical optics, the background for Kummer's work in the 1860s. Such congruences of lines will envelope a caustic surface (\textit {Brennfl\"ache}), which in this case will be a surface of the fourth order and class. Unlike Kummer surfaces, however, Pl\"ucker surfaces have a double line $g$ and other singularities, but never more than 8 nodes (thus Figure \ref{fig:LMS-13} shows the maximal case). 
Pl\"ucker's complex surfaces thus represent a degenerate type of Kummer surface where the nodal line $g$ forms a double line. This double line 
contains four higher point singularities (so-called  
pinch points) where the leaves of the surface meet. These pinch points also demarcate the boundaries between real and imaginary portions of the surface as illustrated in Figure \ref{fig:Pluecker_11112}.

In the posthumously published \cite{Pluecker1869}, which was actually written by Pl\"ucker's student, Felix Klein, one finds a classification scheme for 78 distinct types of so-called equatorial complex surfaces (see below). This scheme is mainly of interest because it shows how Pl\"ucker went about constructing the various cases in a systematic fashion. 
Pl\"ucker also designed prototypes for models of certain types of complex surfaces, both equatorial and meridianal, which were then manufactured by a company in Cologne.

For many of his models,  Pl\"ucker 
 chose the  line 
 $g\subset  E_{\infty}$, the plane at infinity. The family of planes $E(g)$ that contain such a nodal line $g$ are then parallel,  and the quartic $S_4$ enveloped by the lines in $K_2  \cap K_1(g)$ is then an \textit {equatorial surface} (examples are shown in Figures \ref{fig:nr23} and \ref{fig:LMS-13}). All other cases are 
known as \textit {meridianal surfaces}. 
Since each $E(g) \cap S_4$ is a quartic curve that contains $g$ as a double line, these  curves break up into a family of conics $C_2$ together with $g$. In the case of equatorial surfaces, Pl\"ucker showed that the equation for the surface in point coordinates has the form $$\frac{y^2}{Ex^2+2Ux+C}+\frac{z^2}{Fx^2-2Rx+B}+1=0.$$ This coordinate system is chosen so that the $x$-axis coincides with the diameter of the complex $K_2$. This ensures that the conics that form the latitudinal curves (Breitenkurven) $E(g) \cap S_4$ have centers that lie on the $x$-axis. Pl\"ucker nexts shows how to relate the lengths of their semi-diameters to two other fundamental conics in the $xy$- and $xz$-planes. The equations for their respective intersections with $S_4$ are then
$$y^2+Ex^2+2Ux+C=0; \qquad z^2+Fx^2-2Rx+B=0.$$ One notes that the zeroes of the quadratics in $x$ in the denominators correspond to two singular lines that lie in each of the two coordinate planes. Ignoring these, Pl\"ucker calls the two remaining conics the \emph{characteristic curves} of the surface. Depending on their reality and relative position, he shows that there are 17 possible cases.

\begin{figure}[ht]
        \centering 
        \includegraphics[width= 0.5 \textwidth]{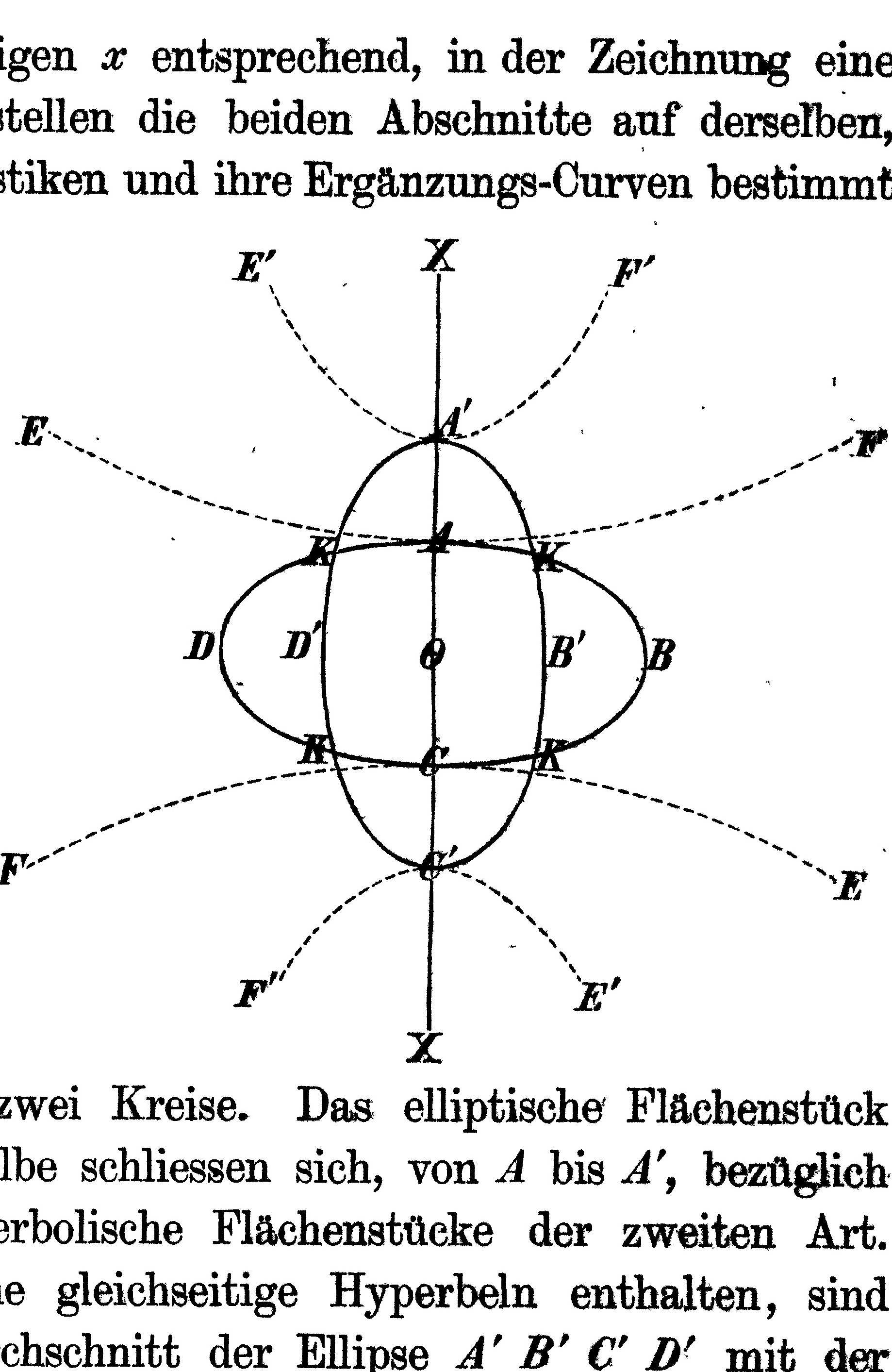}
        \caption{Pl\"ucker's figure of the characteristic curves for an equatorial surface (Pl\"ucker 1868, 351)}
\label{fig:pl-gdr}
\end{figure}

\begin{figure}[ht]
        \centering 
        \includegraphics[width= 0.33 \textwidth]{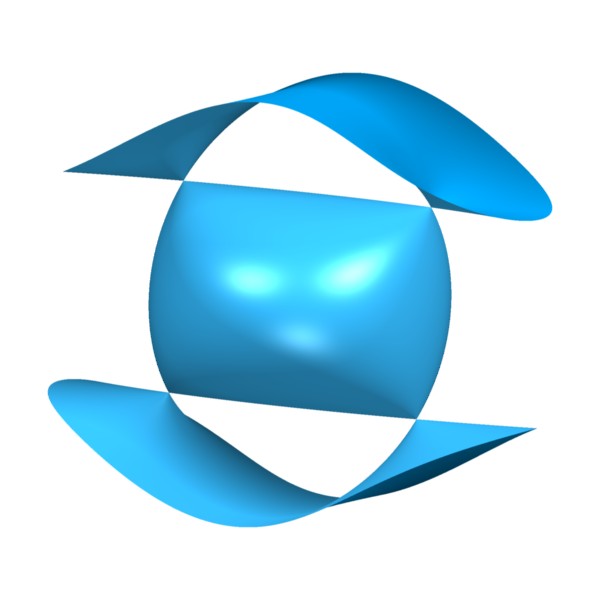}
        \caption{An equatorial surface with 4 real singular points. Graphic by Oliver Labs.}
\label{Aequatorial_Goet_009_110}
\end{figure}

To illustrate how these characteristic curves can be used to construct the latitudinal curves of $S_4$, he discusses the case shown in Figure \ref{fig:pl-gdr}. A graphic of the resulting surface is given in Figure \ref{Aequatorial_Goet_009_110}. (Original models exist in the collections in G\"ottingen (Nr. 110) and in Munich (Nr. 9).) To obtain the plane figure, Pl\"ucker rotates one of the two planes around the $x$-axis so that it coincides with the other. He thus obtains two coplanar conics with four points (real or imaginary) of intersection with this axis. These correspond to the places where the latitudinal curves transition from one type of conic to a different type. The dotted hyperbolas indicate these transitions, so the conics $E(g) \cap S_4$ that lie between $A$ and $C$  will be ellipses, those between $A$ and $A'$ and $C$ and $C'$ are hyperbolas, and those above $A'$ and below $C'$ are imaginary hyperbolas. The points of intersection of the characteristic curves marked $K$ correspond to two circular cross sections, whereas the unmarked intersection points of the inner hyperbola with the elongated ellipse yield cross sections which are rectangular hyperbolas.

Pl\"ucker's entire classification involved 78 cases in all, but the methodology he employed is of far more interest and can be sketched as follows. For each of his 78 cases, Pl\"ucker introduced a special symbol based on his construction method. The case just considered (Nr. 9) has the Pl\"ucker symbol $I_1H_2E_1H_2I_1$, where $E$, $H$, and $I$ indicate the alternation between elliptic, hyperbolic, and imaginary latitudinal curves. The suffix $1$ signifies that the bounding singular lines are parallel, $2$ that they lie in perpendicular directions. If the two characteristic curves are real ellipses, then two other  neighboring cases arise. Nrs. 8 and 10  
occur when the intersection points of the ellipses with the $x$-axis lie outside each other or,  respectively, the points of intersection alternate between the two curves. In the latter case 10, the Pl\"ucker symbol is $I_2H_2E_2H_2I_2$, which differs only slightly from case 9 in that the suffixes are identically 2. So imagining that we view the cross sections kinematically by passing from one type of curve to the next, then each singular line encountered lies in the opposite direction from the one preceding.

From these 17 canonical cases, Pl\"ucker goes on to derive systematically all the other degenerate types. He starts by 
considering cases where two of the singular lines coincide. For example, case Nr. 9 with symbol $I_1H_2E_1H_2I_1$ can pass over into  surface Nr. 23 with symbol $I_1H_2H_2|I_1$, where $|$ denotes a double singular line. Here the latitudinal conics within the visible segment  of the $x$-axis are hyperbolas, where at one end of the segment the double singular line appears. As before, the curves above and below the segment are imaginary. Another type of degeneracy occurs when the two characteristic curves have coincident intersection points on the $x$-axis.
An instance of this is Nr. 34 with Pl\"ucker symbol $I\times E H_2I$, where the symbol $\times$ indicates that the latitudinal curve consists of two perpendicular lines. This case can arise by continuous deformation starting with either of the two cases 9 or 10.

\subsection*{On Deforming Quartics}

Soon after his teacher's death, Klein supplemented Pl\"ucker's collection with four additional models to illustrate the connection between Kummer surfaces and three principal types of Pl\"ucker surfaces. 
Beginning with a Kummer surface, Klein produced these models by considering the relation of the double line $g$ to the quadratic line complex. For $g\notin K_2$, the complex surface has eight double points and four pinch points, whereas for $g\in K_2$ there will be only four double points. Finally, if $g$ is a singular line in $K_2$, then the Pl\"ucker surface will have just two double points.    
 This eventually led to a  program for classifying Pl\"ucker surfaces by means of deformation techniques, an important part of  research on real quartics later pursued by Klein's student Karl Rohn (1855--1920) \cite{Rohn1879}, \cite{Rohn1881}. It seemed evident to Klein that one should start with a  Kummer surface, carrying out  deformations that systematically lowered the number of double points. Thus, the double line $g$ on a Pl\"ucker surface in effect absorbed eight of the 16 double points on a Kummer surface. Only  eight then remain as the others  pass over into the four pinch points on $g$. 
Klein's first model (Figure \ref{fig:Klein-Kummer})  proved very useful for picturing the configuration of  singularities on a Kummer surface. Its connection with the model shown in Figure \ref{fig:LMS-13}, on the other hand, was much less clear. 
One of the main difficulties  in visualizing how a Kummer surface can, via deformation, acquire  a double line is that the Pl\"ucker quartics also contain four other lines as tangential singularities. So in the process of absorbing eight of the sixteen double points to form four pinch points on the double line $g$, four lines  determined by the four tropes that pass through $g$,  need to come into view.

\begin{figure}[h]
        \centering 
        \includegraphics[width=6.0cm,height=7.0cm]{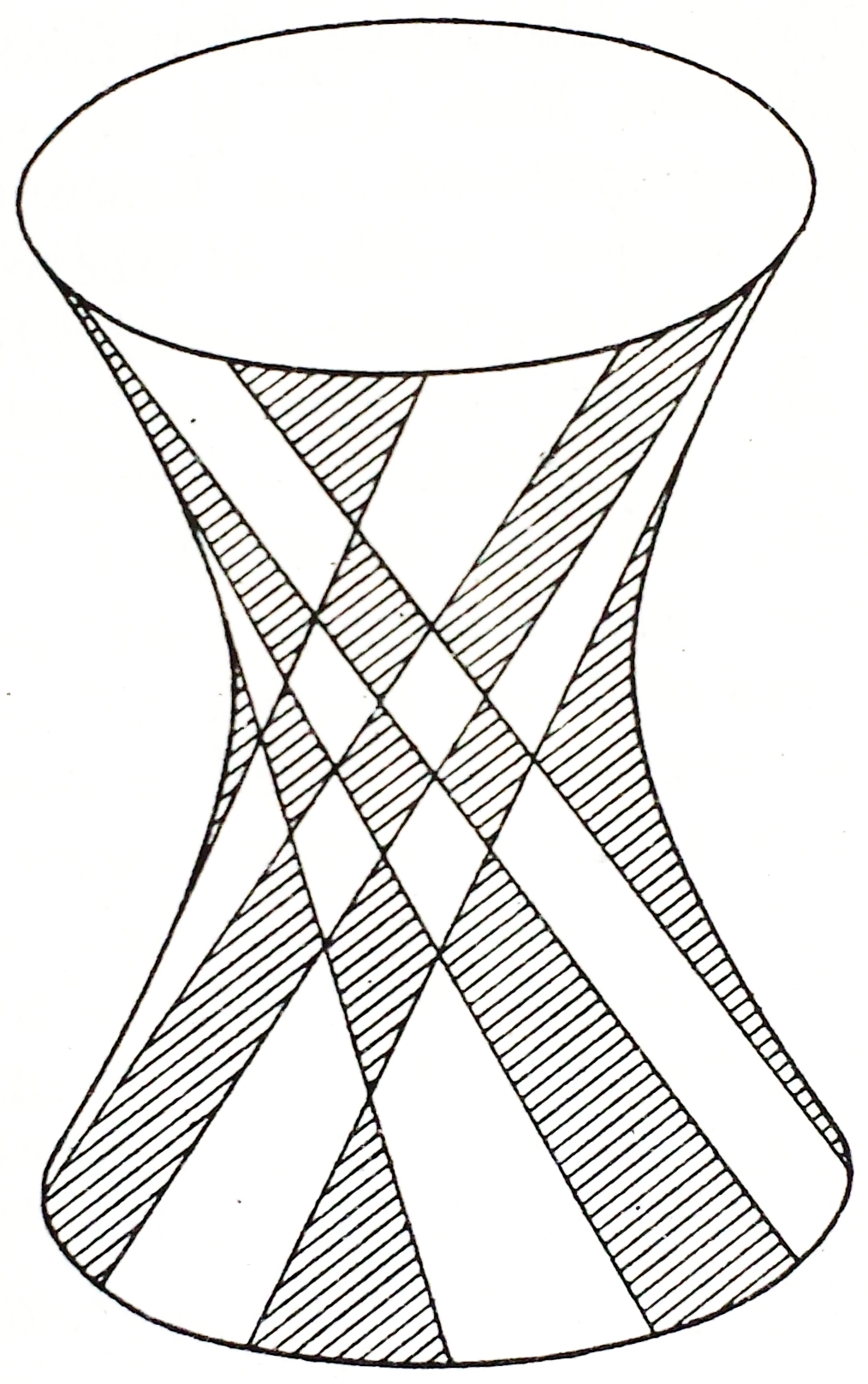}
        \caption{Rohn's Degenerate Kummer Surface \cite{Klein1921-23}, 2: 9}
\label{fig:DoubleHyperboloid_Kummer}
\end{figure}

In the late 1870s, Rohn found a way to link the Pl\"ucker and Kummer surfaces by deformations.
Instead of 
starting with a general Kummer surface, he began with a highly degenerate case of a quartic surface, namely a quadric counted twice (see Figure \ref{fig:DoubleHyperboloid_Kummer}). A familiar quadric surface $F_2$ is a hyperboloid of one sheet, which contains two systems of real generators, or rulings by lines. These one-parameter sets of lines arise naturally in line geometry simply as the intersections of three first-degree line complexes. 

Rohn next took four lines each from the two rulings, which then produces the checkerboard configuration of lines on the surface of the hyperboloid shown in Figure \ref{fig:DoubleHyperboloid_Kummer}. This pattern of lines yields 16 double points, and the shaded portion can then be treated as the region corresponding to the real points on a degenerate Kummer surface, each piece of which
corresponds to a flattened image of the tetrahedra in Rohn's model, discussed above. Passing through infinity we then see that there are exactly eight pieces, in agreement with what we observed with the model in Figure \ref{fig:Rohn_Kummer}. So Rohn could now obtain a general Kummer surface by blowing up the flattened tetrahedra while leaving their vertices fixed.

\begin{figure}[h]
        \centering 
        \includegraphics[width=6.0cm,height=7.0cm]{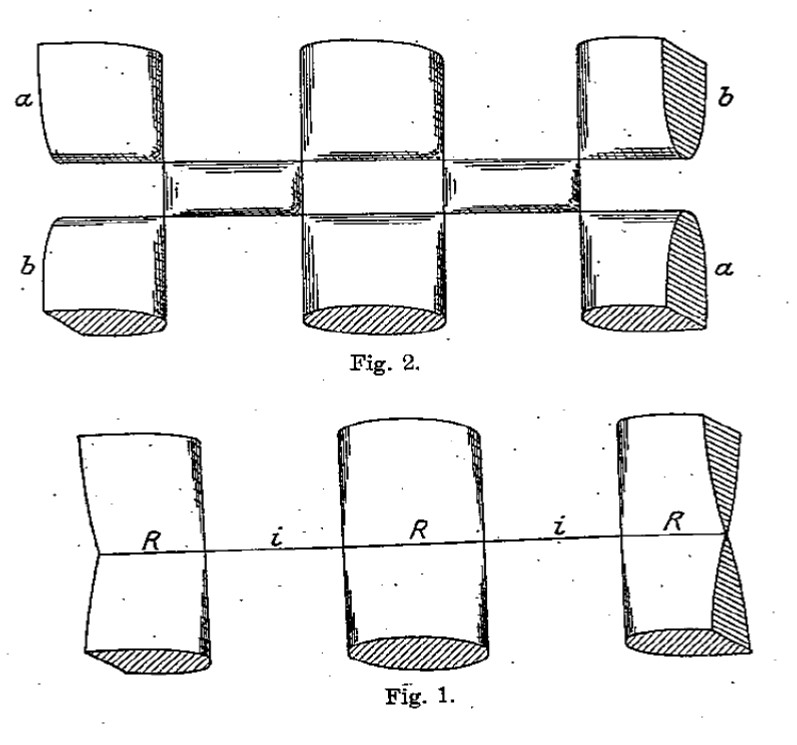}
        \caption{A Kummer surface with 8 nodes passing into a Pl\"ucker surface \cite{Klein1921-23}, 2: 8--9}
\label{fig:Ku-Pl}
\end{figure}

Rohn's approach  made it far easier to visualize the passage from a Kummer surface to a Pl\"ucker quartic with a double line. One can carry out this deformation by starting with the degenerate Kummer surface, then what happens is that two of the four lines in, say, the first ruling gradually fall together. As they do so, the eight singular points on these two lines combine in pairs to form four pinch points on the double line. After this, one only has to blow up the flattened pieces, which are now five instead of eight in number. So by this simple means Rohn was able to show how these degenerate quartics made the relationship between the Kummer and Pl\"ucker surfaces more transparent. This insight  made it possible to visualize 
the relationship between the singularities on these two types of special quartic surfaces. Klein later added a schematic drawing to illustrate how the double line with four pinch points is formed when a Kummer surface passes over into a Pl\"ucker surface (Figure \ref{fig:Ku-Pl}).

The Kummer surfaces -- as well as those quartics related to them but with various other singularities --
turned out to be of central importance for the classification of quadratic line complexes, the topic of Klein's doctoral dissertation. Klein only dealt with the generic case, however, whereas a detailed classification requires using Weierstrass' theory of elementary divisors to analyze all possible degeneracies among the eigenvalues of a 6 x 6 matrix. 
Five years later, Klein's student Adolf Weiler gave the first detailed analysis of 48 different types of quadratic line complexes by making use of their singularity surfaces \cite{Weiler1874}. In the most general case these are Kummer surfaces, which then pass over to Pl\"ucker quartics if the complex contains double lines. In several cases, Weiler found other types of quartic surfaces that had been studied earlier by Luigi Cremona \cite{Cremona1868}, Arthur Cayley, Jakob Steiner, Sophus Lie, and Ludwig Schl\"afli. Thus, his classification scheme drew on much recent knowledge from algebraic surface theory. Ten years later a third doctoral dissertation, written by Corrado Segre, presented still another even more refined classification of quadratic line complexes; Segre's work has remained the last word on this topic \cite{Rowe2017}.

\end{document}